\documentclass{article}
\usepackage{amsmath}
\usepackage{amssymb}
\usepackage{amsfonts}
\usepackage{bm}
\usepackage{float}
\usepackage{graphicx}
\usepackage{hyperref}

\begin{document}

\begin{center}
{\large\bf Mathematical Modeling of Networks in Strength of Materials}

\vspace{0.2in}

Ioannis Dassios$^{1*}$

\vspace{2mm}
\footnotesize
$^1$Aristotle University of Thessaloniki, Greece

\vspace{2mm}
$^*$A Problem in Numerical Analysis for Undergraduate Students
\end{center}

{\footnotesize
\noindent
\textbf{Abstract:} 
We study a material modeled as a network of nodes connected by edges. Using a discrete approach, we build a nonlinear algebraic system that connects applied forces to internal forces and node positions. The model can describe elasticity, plasticity, and possibly cracking. The goal is to solve this system and understand how the material responds. Students are asked to start with a simple triangle example and then apply the method to larger structures. The final aim is to solve the full system and justify the results, leading to a possible publication.
}
\\\\
\noindent{\footnotesize
More info: \href{https://scholar.google.com/citations?user=t4zvXzQAAAAJ&hl=en}{Google Scholar} \\
Personal page: \href{https://ioannisdassios.wordpress.com/publications/}{https://ioannisdassios.wordpress.com/publications/}}

\section{Initial Configuration of the Material}

We consider a material modeled as a network of discrete points (called \emph{nodes}) connected by bonds (called \emph{edges}). This forms a lattice structure.
\\\\
Let \( D \in \mathbb{R}^{n \times 3} \) be the matrix representing the initial positions of the nodes in three-dimensional space. Each row \( D_i \in \mathbb{R}^3 \) represents the position of node \( i \), containing its spatial coordinates along the three directions: length, width, and height.
\\\\
Let \( b \in \mathbb{R}^{m \times 3} \) be the matrix that contains the initial position vectors of the \( m \) bonds in the material. Each row \( b_i \in \mathbb{R}^3 \) corresponds to bond \( i \), and gives the vector from its starting node to its ending node - that is, the direction and length of the bond in the undeformed configuration.
\\\\
These vectors are computed from the nodal positions using the linear system:
\begin{equation}\label{eq1}
b = AD,
\end{equation}
where \( A \in \mathbb{R}^{m \times n} \) is the \textbf{incidence matrix}, which encodes the structure of the graph - specifically, how nodes are connected to form bonds.
\\\\
The incidence matrix is defined as:
\[
a_{ij} = 
\begin{cases}
1, & \text{if node } j \text{ is the start of bond } i, \\
-1, & \text{if node } j \text{ is the end of bond } i, \\
0, & \text{otherwise}.
\end{cases}
\]
Each row of \( A \) contains exactly one \( +1 \) and one \( -1 \), corresponding to the two nodes connected by each bond. In this way, the equation \( b = AD \) captures both the topology and the geometry of the material in its undeformed (reference) configuration, formulated as a system of linear equations.
\\\\
To summarize:
\begin{itemize}
  \item \( D \): initial positions of the nodes (known - they describe the undeformed material),
  \item \( b \): initial position vectors of the bonds (also known),
  \item \( A \): incidence matrix (connectivity of the graph, to be determined),
  \item \( b = AD \): linear system that relates node positions and bond vectors,
\end{itemize}

\section{Unknown Configuration Due to Applied Forces}

Now we assume that an external force is applied to a specific part of the material. As a result, the configuration of the material changes: the nodes may shift position in response to the load, depending on the physical behavior of the bonds. This response may include elastic or plastic deformation, or even bond breakage, depending on the force magnitude and the material's properties.
\\\\
Let \( X \in \mathbb{R}^{n \times 3} \) be the matrix of \textbf{updated positions} of the nodes. Each row \( X_i \in \mathbb{R}^3 \) gives the new coordinates of node \( i \) after the application of external forces. These updated positions may result from different types of material response, including elastic or plastic behavior, or in some cases, even bond failure. 
\\\\
Let \( y \in \mathbb{R}^{m \times 3} \) be the matrix of \textbf{updated position vectors of the bonds}. Each row \( y_i \in \mathbb{R}^3 \) corresponds to bond \( i \), and represents its new (unknown) position vector after the material responds to the applied forces.
\\\\
These updated bond vectors are related to the updated nodal positions \( X \in \mathbb{R}^{n \times 3} \) through the linear system:
\begin{equation}\label{eq2}
AX=y,
\end{equation}
where \( A \in \mathbb{R}^{m \times n} \) is the incidence matrix that encodes the connectivity of the graph. This equation expresses the fact that each bond vector is determined by the difference between the positions of the two nodes it connects, as described by the structure of \( A \).
\\\\
In our setup, a force is applied only to part of the material. As a result, the positions of some nodes remain unchanged because they are unaffected by the applied force. These positions are known. The remaining nodes respond to the force, and their positions are unknown and need to be computed. To organize this, we split the nodal positions into two groups: the unknown ones and the known ones.
\\\\
The matrix \( X \in \mathbb{R}^{n \times 3} \) contains the updated positions of all \( n \) nodes. It can be written as:
\[
X =
\begin{bmatrix}
X_1 \\
X_2 \\
\vdots \\
X_p \\
X_{p+1} \\
X_{p+2} \\
\vdots \\
X_n
\end{bmatrix},
\]
where each row \( X_i \in \mathbb{R}^3 \) gives the position of node \( i \) in the updated configuration.
\\\\
We now distinguish between:
\begin{itemize}
  \item the first \( p \) nodal positions, \( X_1, X_2, \dots, X_p \), which are \textbf{unknown}, and
  \item the remaining \( q = n - p \) positions, \( X_{p+1}, \dots, X_n \), which are \textbf{known} due to applied boundary conditions or controlled displacements.
\end{itemize}
For compactness, we group the unknown and known nodal positions into the block vectors:
\[
X = 
\begin{bmatrix}
\mathbf{X}_P \\
\mathbf{X}_Q
\end{bmatrix},
\]
where:
\[
\mathbf{X}_P =
\begin{bmatrix}
X_1 \\
X_2 \\
\vdots \\
X_p
\end{bmatrix}
\in \mathbb{R}^{p \times 3}, 
\qquad
\mathbf{X}_Q =
\begin{bmatrix}
X_{p+1} \\
X_{p+2} \\
\vdots \\
X_n
\end{bmatrix}
\in \mathbb{R}^{q \times 3}.
\]
where:
\begin{itemize}
  \item \( X_1, \dots, X_p \in \mathbb{R}^3 \) are the \textbf{unknown} nodal positions (the free part of the material), and are collectively denoted by \( \mathbf{X}_P \in \mathbb{R}^{p \times 3} \),
  \item \( X_{p+1}, \dots, X_n \in \mathbb{R}^3 \) are the \textbf{known} nodal positions (where forces are directly applied), and are collectively denoted by \( \mathbf{X}_Q \in \mathbb{R}^{q \times 3} \), with \( q = n - p \).
\end{itemize}
The goal is to determine the unknown nodal positions \( \mathbf{X}_P \).

\section{Forces and Their Relation to Deformation}

The matrix \( F \in \mathbb{R}^{m \times 3} \) contains the forces transmitted through the bonds. Each row \( F_i \in \mathbb{R}^3 \) represents the force acting along bond \( i \).
\\\\
Similarly, recall that \( y \in \mathbb{R}^{m \times 3} \) contains the deformed bond vectors:
\[
y = 
\begin{bmatrix}
y_1 \\
y_2 \\
\vdots \\
y_m
\end{bmatrix}, \qquad
F = 
\begin{bmatrix}
F_1 \\
F_2 \\
\vdots \\
F_m
\end{bmatrix}.
\]
We assume that the force in each bond acts along its deformed direction. That is, the force vector \( F_i \) is parallel to the bond vector \( y_i \). This means:
\[
F_i = \frac{|F_i|}{|y_i|} y_i,
\]
where \( |F_i| \) is the scalar magnitude of the internal force in bond \( i \), and \( |y_i| \) is the length of the deformed bond vector.
\\\\
Equivalently, this relationship can be written as:
\begin{equation}\label{eq3}
F = K(y) y, 
\end{equation}
where \( K(y) \in \mathbb{R}^{m \times m} \) is a diagonal matrix that depends on the deformation \( y \), with entries:
\[
K(y) = \mathrm{diag}\left( \frac{|F_1|}{|y_1|}, \frac{|F_2|}{|y_2|}, \dots, \frac{|F_m|}{|y_m|} \right).
\]
\begin{figure}[H]
    \centering
    \includegraphics[width=0.6\textwidth]{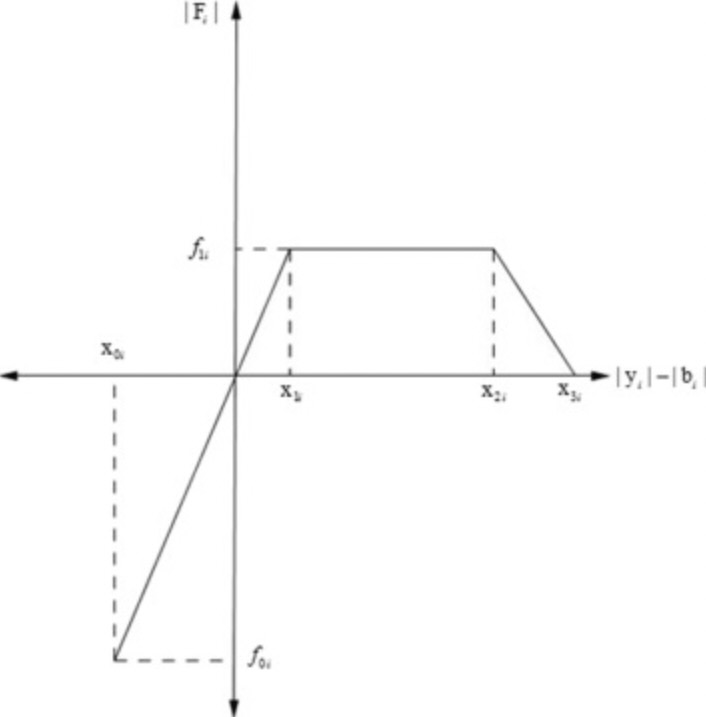}
    \caption{Graph showing the relationship between the force magnitude \( |F_i| \) and the extension \( |y_i| - |b_i| \) of a bond.}
    \label{Figure 1.jpeg}
\end{figure}
The magnitude \( |F_i| \) is determined by a material law that relates the force to the change in bond length, such as:
\[
|F_i| = f(|y_i| - |b_i|),
\]
where \( |b_i| \) is the initial (undeformed) length of bond \( i \), and \( f \) is a given 
function that models the elastic or plastic response of the material. This relationship captures how each bond resists stretching or compression based on how much its length has changed from the reference configuration.
\\\\
Using experimental data or material testing in the laboratory, we can determine the relationship between the force magnitude \( |F_i| \) and the change in bond length \( |y_i| - |b_i| \). This relation depends on the physical properties of the material and may include elastic, plastic, or damage behavior. The typical shape of this function is illustrated in the graph below (see Figure~\ref{Figure 1.jpeg}).

\section{Reaction Forces}

When we apply a force to part of the material, that part responds by generating internal forces at the nodes. These are called \textbf{reaction forces}. Let \( B \in \mathbb{R}^{n \times 3} \) be the matrix of all reaction forces:
\[
B =
\begin{bmatrix}
B_1 \\
B_2 \\
\vdots \\
B_p \\
B_{p+1} \\
\vdots \\
B_n
\end{bmatrix},
\]
where each \( B_i \in \mathbb{R}^3 \) is the total force acting on node \( i \). We divide the nodes into two groups:
\begin{itemize}
  \item Nodes \( 1 \) to \( p \) are the ones where we directly apply external forces. So the corresponding reaction forces \( B_1, \dots, B_p \) are \textbf{known}.
  \item Nodes \( p+1 \) to \( n \) are free to move and are not directly loaded. The corresponding reaction forces \( B_{p+1}, \dots, B_n \) are \textbf{unknown}, and we have to compute them.
\end{itemize}

\medskip

We now group these into two blocks:
\[
B =
\begin{bmatrix}
\mathbf{B}_P \\
\mathbf{B}_Q
\end{bmatrix},
\quad \text{where} \quad
\mathbf{B}_P =
\begin{bmatrix}
B_1 \\
\vdots \\
B_p
\end{bmatrix},
\quad
\mathbf{B}_Q =
\begin{bmatrix}
B_{p+1} \\
\vdots \\
B_n
\end{bmatrix}.
\]

Here, \( \mathbf{B}_P \in \mathbb{R}^{p \times 3} \) is \textbf{known}, and \( \mathbf{B}_Q \in \mathbb{R}^{q \times 3} \) is \textbf{unknown}, with \( q = n - p \).
The final step is to write the equation that connects the forces \( F \), which act inside the material (along bonds), to the external reaction forces \( B \), acting at the nodes. This relationship is given by the equation:
\begin{equation}\label{eq4}
A^T F = B.
\end{equation}
In this equation, \( A^T \) is the transpose of the incidence matrix, which tells us how the bonds are connected to the nodes. The matrix \( F \) contains the internal forces acting along the bonds, while \( B \) contains the total external force applied at each node.

\section{What Needs to Be Solved}

We now show how the final equation is formed by combining the previous steps. Substituting \eqref{eq2} into equation \eqref{eq3}, we get:
\[
F = K(y)\, y = K(AX)\, AX.
\]
Substituting \eqref{eq4} for \( F \), we obtain:
\begin{equation}\label{eq5}
A^T K(AX)\, AX = B.
\end{equation}
This is the final nonlinear system that links the geometry of the structure, the internal forces in the material, and the external forces applied at the nodes.
The unknowns in this equation are the unknown positions collected in \( \mathbf{X}_P \), and the unknown reaction forces at those nodes, collected in \( \mathbf{B}_Q \). 
\\\\
We are given the initial configuration of the material: the positions of the nodes \( D \), and the initial bond vectors \( b \). From this data, we compute the incidence matrix \( A \), which encodes how the nodes are connected.
\\\\
The main system we need to solve is \eqref{eq5}:
\[
A^T K(AX)\, AX = B.
\]
We now recall how the matrices \( X \) and \( B \) are split into known and unknown parts:
\[
X = 
\begin{bmatrix}
\mathbf{X}_P \\
\mathbf{X}_Q
\end{bmatrix}, \qquad
B = 
\begin{bmatrix}
\mathbf{B}_P \\
\mathbf{B}_Q
\end{bmatrix}.
\]
The matrix \( \mathbf{X}_Q \) is known, and \( \mathbf{B}_Q \) is also known.  
The matrix \( \mathbf{X}_P \) is unknown, and \( \mathbf{B}_P \) is also unknown.

The matrix \( K(AX) \) is diagonal and is based on the force-extension graph in Figure~\ref{Figure 1.jpeg}.
\\\\
\textbf{Unknowns of the problem:}
\[
\boxed{\mathbf{X}_P \quad \text{and} \quad \mathbf{B}_P}
\]
We solve equation \eqref{eq5} using the known \( \mathbf{X}_Q \) and \( \mathbf{B}_Q \).

\section*{Conclusion and What Comes Next}

To be able to publish these results in an international, peer-reviewed journal of high impact, two main steps are needed:
\\\\
First, we must solve the nonlinear system given in equation~\eqref{eq5}. This can be done using the numerical method we introduced for solving nonlinear algebraic equations. That method can also be expanded or generalized to handle nonlinear algebraic systems like the one we have here.
\\\\
Second, we need to support and explain the results with a numerical example. A good choice is to use a 3D structure like an octahedron with 15 nodes and 14 edges. This example is large enough to show meaningful internal forces and demonstrate how the method works in a realistic case.
\\\\
Before solving the full system in equation~\eqref{eq5}, it's helpful to start with a small example: 3 nodes and 3 edges forming a triangle in 3D. For example, take the initial positions:
\[
D_1 = (1, 1, 0), \quad D_2 = (2, 1, 1), \quad D_3 = (1, 2, 1).
\]
Connect the nodes with 3 edges: between nodes \( D_1 \)-\( D_2 \), \( D_2 \)-\( D_3 \), and \( D_3 \)-\( D_1 \). Construct the incidence matrix \( A \), and apply external forces to two of the nodes - for example, to \( D_1 \) and \( D_2 \) - and fix the position of \( D_3 \). Use the force-extension graph shown in Figure~\ref{Figure 1.jpeg} to define the internal force in each bond, based on how much it stretches.Then use the equations \eqref{eq2}, \eqref{eq3}, \eqref{eq4}, to build the nonlinear system \eqref{eq5}. This leads to a simple system of equations where the unknowns are the new positions of nodes \( X_1 \) and \( X_2 \). Once you solve this small case, you can generalize the method to handle larger, more complex networks, and the same idea can be extended to solve the full system~\eqref{eq5}.
\\\\
If anyone would like to try out an example, discuss the method further, or ask questions, feel free to contact me.

\end{document}